\documentclass[12pt]{article}
\usepackage[textwidth=15cm, textheight=22cm]{geometry}

\usepackage{graphicx}  
\usepackage{dcolumn}   
\usepackage{bm}        

\usepackage{amsfonts,epsfig,amsmath,amsthm,amssymb,verbatim}
\usepackage{color}

\usepackage{soul}
\usepackage{authblk}
\usepackage[onehalfspacing]{setspace}
\usepackage{nicefrac}
\usepackage{todonotes}
\usepackage{subfig}
\usepackage{epstopdf}
\usepackage{hyperref}

\def\R{\mathbb{R}}

\def\cC{\mathcal{C}}

\def\cH{\mathcal{H}}

\def\cO{\mathcal{O}}

\def\cS{\mathcal{S}}

\def\cV{\mathcal{V}}

\def\txta{{\textnormal{a}}}

\def\txtd{{\textnormal{d}}}
\def\txte{{\textnormal{e}}}

\def\txtr{{\textnormal{r}}}

\def\ra{\rightarrow}
\def\I{\infty}

\newcommand{\be}{\begin{equation}}
\newcommand{\ee}{\end{equation}}
\newcommand{\benn}{\begin{equation*}}
\newcommand{\eenn}{\end{equation*}}
\newcommand{\bea}{\begin{eqnarray}}
\newcommand{\eea}{\end{eqnarray}}
\newcommand{\beann}{\begin{eqnarray*}}
\newcommand{\eeann}{\end{eqnarray*}}

\DeclareMathOperator{\var}{Var}

\begin{document}

\title{Warning signs for non-Markovian bifurcations: colour blindness and scaling laws}
\author[1,2]{Christian Kuehn}
\author[1]{Kerstin Lux}
\author[3]{Alexandra Neam\c tu}
\affil[1]{Technical University of Munich, Department of Mathematics, Boltzmannstr.~3, 85748 Garching b.~M\"unchen, Germany}
\affil[2]{Complexity Science Hub Vienna, Josefst\"adter Str.~39, 1080 Vienna, Austria}
\affil[3]{University of Konstanz, Department of Mathematics and Statistics,  Universit\"atsstr.~10, 78464 Konstanz, Germany}

\maketitle

\begin{abstract}
Warning signs for tipping points (or critical transitions) have been very actively studied. Although the theory has been applied successfully in models and in experiments for many complex systems {such as} for tipping in climate systems, there are ongoing debates, when warning signs can be extracted from data. In this work, we {shed light on this debate by considering different types of underlying noise. Thereby, we }
significantly advance the general theory of warning signs for nonlinear stochastic dynamics. A key scenario deals with stochastic systems approaching a bifurcation point dynamically upon slow parameter variation. The stochastic fluctuations are generically able to probe the dynamics near a deterministic attractor to {reveal} critical slowing down. Using scaling laws {near bifurcations}, one can then anticipate the distance to a bifurcation. Previous warning signs results assume that the noise is Markovian, most often even white. Here, we study warning signs for non-Markovian systems including coloured noise and $\alpha$-regular Volterra processes (of which fractional Brownian motion and the Rosenblatt process are special cases). We prove that early warning scaling laws can disappear completely or drastically change their exponent based upon the parameters controlling the noise process. This provides a clear explanation, why applying standard warning signs results to reduced models of complex systems may not agree with data-driven studies. We demonstrate our results numerically in the context of a box model of the Atlantic Meridional Overturning Circulation (AMOC). 
\end{abstract}


\textbf{Introduction.} Predicting drastic, sudden, and large changes in dynamics, so-called tipping points or critical transitions, is a major challenge across all sciences~\cite{Helbing,Schefferetal,Schefferetal2}. Recent results have successfully employed warning signs in cases ranging from small {scale experiments}~\cite{DrakeGriffen} to {large scale ecosystems}~\cite{Carpenteretal1}. The most common framework to extract warning signs is to exploit critical slowing down near bifurcations~\cite{KuehnCT1,Wiesenfeld1}, i.e., a reduced return rate to a steady state (or a more complicated attractor) upon approaching the bifurcation point under parameter variation. Although it is sometimes possible to detect critical slowing down {through perturbation} in experiments~\cite{LimEpureanu}, in many cases this is {infeasible}, e.g., in climate dynamics~\cite{Lenton}, social networks~\cite{KuehnMartensRomero} or for ecosystems~\cite{GuttalJayaprakash}. Yet, {in the presence of noise}, one indirectly detects critical slowing down as growing variance and/or {autocorrelation}~\cite{CarpenterBrock,Wiesenfeld1}. There is a well-developed mathematical theory supporting this approach~\cite{AshwinWieczorekVitoloCox,KuehnCT2}. However, for many complex systems one finds debates, whether warning signs can be extracted from certain observables. A representative example are {tipping points in the climate system}, where arguments for and against the existence of warning signs in time series data have been put forward~\cite{Dakosetal,DitlevsenJohnsen}. Similar discussions can be found in the context of complex ecological systems~\cite{BoettingerRossHastings}.

It is crucial to note that previous arguments relied upon several assumptions. First, one assumes that the complex system can be reduced, i.e., one postulates the existence of a closed dynamical system, where remaining degrees of freedom are addressed by considering stochastic terms~\cite{Gardiner}. Second, the noise in the system is assumed to be white in time~\cite{GottwaldMelbourne,Justetal,JustKantzRoederbeckHelm}, e.g., allowing for the use of classical fluctuation-dissipation theorems~\cite{DitlevsenJohnsen,Kubo}. From a physical perspective, a crucial point is disregarded. Although many complex systems observables certainly satisfy reduced stochastic systems, it is expected that white noise models can fall short of capturing the remaining degrees of freedom. This is a known insight, e.g., from the abstract Mori-Zwanzig~\cite{ChorinHaldKupferman,Ghil.2021,HijonEspanolVanden-EijndenDelgado-Buscalioni,Mori,Zwanzig} reduction framework as well as concrete techniques to reduce large-scale stochastic dynamical systems~\cite{EvansGrigoliniParravicini,KondrashovChekrounGhil}. In particular, the Mori-Zwanzig formalism relies crucially on projection operators, which can often only capture the full system dynamics if in-time memory integrals are introduced in the reduced system, so that the reduced model becomes non-Markovian. Also in many concrete applications, non-Markovian terms appear more directly. For example, in the AMOC, it is physically evident that small-scale box models for observables do not take into account the effects of long-term historic fluctuations induced by sea ice if one uses white noise~\cite{Dijkstra}. And neither do most models for ecological observables take into account historical climate fluctuations leading to pattern formation in the landscape. {Moreover, in their recent study \cite{vanderBolt.2018}, the authors show that climate memory might vary over time including evidence that in ecosystems such as forests and coral reefs a higher climate memory could come along with an increased risk of critical transitions.}

This makes it evident, why one has to study warning signs for history-dependent non-Markovian processes, which is the main contribution of this work. Since it is often difficult to decide from data, which type of memory is present, we analyse various time-correlated noise processes. We provide a general theory calculating scaling laws near bifurcations for coloured noise and for $\alpha$-regular Volterra processes (of which fractional Brownian motion and the Rosenblatt process are special cases).

{The importance of analyzing early warning signs (EWS) in the presence of coloured noise has for example already been pointed out in \cite{Munch.2012}. Therein, the authors emphasize upon the poor performance of classical early warning indicators in the context of ecological regime shifts in the presence of large correlated fluctuations in extrinsic noise. In \cite{Qin.2018}, the authors find that environmental noise can change and possibly reduce the EWS. Furthermore, in \cite{Dutta.2018}, the authors conclude that the ``EWS performance is highly sensitive to noise type'' \cite[p.\ 1253]{Dutta.2018}. They consider in their study, amongst other, a red-shifted Gaussian noise. In contrast to particular systems in \cite{Dutta.2018,Munch.2012,Qin.2018}, we are interested in broader classes of time-correlated noise processes and we also study and compare the related scaling laws for different time-correlated processes near bifurcations. In particular, our study includes Ornstein-Uhlenbeck noise as well as $\alpha$-regular Volterra processes.}

{The $\alpha$-regular Volterra processes are a very broad and important class of stochastic processes used to depict long range dependence (long memory). They find applications in hydrology, and telecommunications (see \cite{Taqqu.1974}) and in finance. In \cite{Cont.2005} an overview is given of the concept of self-similarity and long range dependence in financial modeling. These ideas also have applications in climate science as well. In \cite{Rypdal.2010} a link between sun-climate complexity and integrated global temperature anomalies is described by a fractional Brownian motion.} {Moreover, in \cite{Caraballo.2018}, a stochastic model from epidemiology with fractional Brownian motion is studied to account for long range memory.}

{The Rosenblatt process provides a powerful way of modeling long-range dependence in the underlying physical processes if the assumption of a Gaussian process is not suitable. In this case, one could use the Rosenblatt process instead of the fractional Brownian motion. However, since its analysis and simulation is by far more involved compared to that of fractional Brownian motion, its use in applied sciences (see e.g.~\cite{Tudor.2009} for a binary market model based upon the Rosenblatt process) seems to be limited so far.}
{In our analysis, w}e find that classical variance scaling laws for white noise change drastically for non-Markovian cases, or even disappear completely. We provide a detailed dependence of the warning signs upon the history-dependent parameters such as correlation time, $\alpha$-regularity, and self-similarity. We demonstrate our results numerically for the Stommel-Cessi model of the AMOC~\cite{Cessi.1994,Stommel1}.

{Discrepancies between data analysis and theoretical predictions can appear for classical EWS~\cite{Dakosetal,DitlevsenJohnsen}. One alternative to address this issue~\cite{Boers.2018} is to use additional indicators, e.g., based on wavelet analysis such as a local Hurst exponent as a suitable replacement of the classical autocorrelation indicator. Besides further developing possibly better suited warning indicators, it is also important to shed light on the reasons for a failure of classical warning indicators. Here, our results on the dependence of the warning signs upon the history-dependent parameters could} provide a natural explanation, why {sometimes no} agreement on {predictability of} events from data via {classical} warning signs is available.\medskip 

\textbf{Background Theory.} We briefly describe the background of warning signs for bifurcation-induced critical transitions. Consider a single macroscopic observable $x=x(t)$ together with a slowly-driven variable $y=y(t)$ within the standard form of a fast-slow system~\cite{KuehnCT1,KuehnCT2}
\be
\label{eq:fs}
\begin{array}{lclcl}
\frac{\txtd x}{\txtd t}  &=:& x' &=& f(x,y,\varepsilon),\\
\frac{\txtd y}{\txtd t}  &=:& y' &=& \varepsilon g(x,y,\varepsilon),
\end{array}
\ee
where $f,g:\R^2\times[0,\varepsilon_0]\ra \R$ for some small $\varepsilon_0>0$ are smooth maps, and $\varepsilon>0$ is a small parameter. For $\varepsilon= 0$ in~\eqref{eq:fs}, one obtains the fast subsystem
\be
\label{eq:fss}
x' = f(x,y,0),\quad y' = 0,
\ee
where $y$ can be viewed as a parameter. The steady states of~\eqref{eq:fss} form the critical manifold $\cC_0:=\{(x,y)\in\R^2:f(x,y,0)=0\}$. A subset $\cS_0\subset \cC_0$ is called normally hyperbolic if $\partial_x f(p,0)\neq 0$ for all $p\in \cS_0$. It is called attracting (resp.~repelling) if $\partial_x f(p,0)<0$ (resp.~$\partial_x f(p,0)>0$) for all $p\in \cS_0$. The critical manifold is also the domain of the slow subsystem
\be
\label{eq:sss}
0 = f(x,y,0),\quad \frac{\txtd y}{\txtd s}  =:\dot{y} = g(x,y,0),
\ee
which is obtained from~\eqref{eq:fs} by a time re-scaling $s=\varepsilon t$ and then considering the singular limit $\varepsilon=0$. A natural way to model and study bifurcation-induced critical transitions is to select various normal forms for the fast subsystem. The classical codimension-one example is the fold with $f_{\textnormal{lp}}(x,y)=x^2+y$. For the fast subsystem, a fold occurs at $(x,y)=(0,0)$, where two steady states $x_*=\pm\sqrt{-y}$ collide and annihilate.

Assuming symmetry or existence of a trivial branch of solutions, other important codimension-one normal forms are the transcritical bifurcation $f_{\textnormal{tc}}(x,y)=x(x+y)$ and the (subcritical) pitchfork bifurcation $f_{\textnormal{pf}}(x,y)=x(y+x^2)$. For the full fast-slow system with $0<\varepsilon\ll1$, the fold, transcritical and subcritical pitchfork can all induce critical transitions upon variation of $y$ through the bifurcation point at $y=0$. For example, consider $f=f_{\textnormal{pf}}$, $g=1$, and take $(x(0),y(0))=(x_0,y_0)$ with $x_0\approx 0, x_0\neq 0$, and $y_0<0$. Then a trajectory gets attracted quickly to an  $\cO(\varepsilon)$-neighbourhood of the normally hyperbolic attracting part of the critical manifold $\cC_{0,\textnormal{pf}}^\txta=\{x=0,y<0\}$. However, the manifold $\cC_{0,\textnormal{pf}}^\txtr=\{x=0,y>0\}$ is normally hyperbolic repelling and there are no attractors of the fast subsystem $x'=x(y+x^2)$ for $y>0$ so that the trajectory makes a large excursion after departing from the region $\{x\approx 0\}$. Similar considerations apply to the fold and transcritical fast-slow bifurcations. These bifurcations exhibit critical slowing down. Considering the attracting manifolds below the bifurcation denoted by $\cC_0^\txta$, local stability decreases as $y\nearrow 0$ since
\be
\partial_x f_{\textnormal{pf}}|_{\cC_0^\txta}=y, \quad 
\partial_x f_{\textnormal{tc}}|_{\cC_0^\txta}=y,\quad 
\partial_x f_{\textnormal{lp}}|_{\cC_0^\txta}=-2\sqrt{-y}. \label{eq:a_pf_tc_lp}
\ee

\textbf{White noise.} The above theory extends to stochastic fast-slow systems. For small noise, $0<\sigma\ll1$, the noise picks up a warning sign of the bifurcation induced by critical slowing down. Indeed, consider the stochastic fast-slow system~\cite{Gentz.2006}
\be
\label{eq:fsnoisy}
x' = f(x,y,\varepsilon) + \sigma \frac{\txtd W}{\txtd t},\quad
y' =\varepsilon,
\ee
where $0<\sigma\ll 1$, $W=W(t)$ is a standard Brownian motion, and $f:\R\times \R\ra \R$ is again a suitable normal form. The attracting part of the critical manifold $\cC_0^\txta$ is given by a graph $\{x=h_0^\txta(y)\}$. By solving the slow dynamics of~\eqref{eq:fsnoisy}, linearizing around a trajectory lying inside the attracting critical manifold, and changing to the slow time scale, we get
\begin{align}\label{ly:bm}
\varepsilon \dot{X}=\partial_x f(h_0^\txta(s),s) X + \sigma \sqrt{\varepsilon}\frac{\txtd W}{\txtd s},
\end{align}
{where $s=\varepsilon t$ is the slow time.}
$X$ is a non-autonomous Ornstein-Uhlenbeck (OU) process. Its variance satisfies the Lyapunov ODE~\cite{Gentz.2006,KuehnCT2}
\benn
\varepsilon \dot{V}=2\partial_x f(h_0^\txta(y),y) V + \sigma^2.
\eenn
This fast-slow system (where the time variable $s$ and the slow variable $y$ in this non-autonomous system can be used interchangeably) has a normally hyperbolic attracting critical manifold given by
\benn
\cV_0=\left\{V=-\frac{\sigma^2}{2\partial_x f(h_0^\txta(y),y)}=:H_0(y)\right\}.
\eenn
The adiabatic time-asymptotic variance $\lim_{t\ra \I}V(t)=:V_\I$ can be calculated to leading-order as $V_\I=H_0(y)$. Since $\partial_x f(h_0^\txta(y),y)$ vanishes at the bifurcation points discussed above, we obtain the variance divergence rates
\be
V_{\I,\textnormal{pf}}(y)=\cO(y^{-1})=V_{\I,\textnormal{tc}}(y),~ 
V_{\I,\textnormal{lp}}(y)=\cO(|y|^{-1/2}), \label{eq:scaling_BM}
\ee
as $y\nearrow 0$ along the attracting part of the critical manifold. Hence, increasing fluctuations provide an EWS. The exponents of the divergence can be viewed as universal critical exponents for the respective fast-slow stochastic normal forms.\medskip

\textbf{Non-Markovian/Non-White Noise.} Next, we show that warning sign theory has highly non-trivial practical extensions to time-correlated noise processes. This theory is of paramount importance as warning signs can be misinterpreted. A natural starting point is to consider small additive noise. Indeed, as we are interested in local transition points, generically within the class of smooth diffusion terms, these terms do not vanish at a point, so that a local approximation by small additive noise is mathematically rigorously justified~\cite{Gentz.2006,KuehnCT2}. Under additional assumptions in applications, e.g., the existence of a trivial zero solution, multiplicative noise would be interesting~\cite{LiSieber}. Hence, we focus on time-correlated additive noise. We start with {Ornstein-Uhlenbeck processes $B=\left(B_t\right)_{t\geq 0}$ and \emph{coloured noise} $B'=\left(B'_t\right)_{t\geq 0}$~\cite{Luczak}.} Consider the fast-slow SDEs~\cite{BerglundGentz5}
\be
\label{eq:fsnoisycor1}
\begin{array}{lcl}
x' &=& f(x,y,\varepsilon) +\sigma B_t',\\
B' &=& -\frac{1}{\tau}B + W_t',\\
y' &=& \varepsilon,
\end{array}
\ee
where $\tau>0$ {corresponds to the strength of time correlation in the noise. Note that, in the limit $\tau\rightarrow\infty$, we are back in the white noise case and for $\tau\rightarrow 0$, we have strongly correlated noise with an immediate mean reversion in the limit resembling deterministic dynamics. Here, we fix $\tau$ and analyse the effect of time correlation in the coloured noise.}
Using linearization 
along an attracting critical manifold we get~\cite{BerglundGentz5}
\be
\label{eq:fsnoisycor1a_vec}
\begin{array}{lcl}
	\begin{pmatrix}
		\varepsilon \dot{\xi}\\
		\dot{B}
	\end{pmatrix} &=&
	\underbrace{\begin{pmatrix}
		a(y) & \frac{-\sigma}{\tau}\\
		0 & \frac{-1}{\tau}
	\end{pmatrix}}_{=A}\begin{pmatrix}
	\xi\\
	B
\end{pmatrix}
+\underbrace{\begin{pmatrix}
	\sigma\\
	1
\end{pmatrix}}_{=D}\dot{W}_s.
\end{array}
\ee
where $a(y)=\partial_x f(h_0^\txta(y),y)$ {and $\xi_t = x(t) - h_\varepsilon^\txta(y(t))$, where $\{x=h_\varepsilon^\txta(y)=h_0^\txta(y)+\mathcal{O}(\varepsilon)\}$ is the attracting slow manifold.} The covariance $\textnormal{Cov}_\I$ solves to leading-order the matrix Lyapunov equation $0=A \textnormal{Cov}_\I + \textnormal{Cov}_\I A^\top + DD^\top${; see~\cite[(5.1.11)]{BerglundGentz} and \cite{Bellman}}. Direct algebra yields
\benn
\textnormal{Cov}_\I=\left(\begin{array}{cc}\frac{\sigma^2}{2(1/\tau+|a(y)|)} & 
\frac{\sigma}{2(1/\tau+|a(y)|)} \\ \frac{\sigma}{2(1/\tau+|a(y)|)} & \tau/2 
\end{array}\right).
\eenn
The variance component of the main observable $x$ gives
\be
\lim_{y\nearrow 0 }V_\I=\frac{\sigma^2}{2(1/\tau+|a(0)|)}=\frac{\sigma^2 \tau}{2} \label{eq:scalingLaw_color}.
\ee
{As we consider the singular limit case $\varepsilon=0$ with a small enough noise intensity $\sigma$ and fixed $\tau$, the variance does not diverge but increases to the bound $\sigma^2\tau/2$,} as long as $\tau$ does not dominate $\lim_{y\nearrow 0}\lvert a(y)\rvert$, i.e, as long as $\tau$ is fixed. Hence, noise modelled by a time-correlated process can make us \emph{``colour blind''} {in the sense that we might not infer the distance to a bifurcation point based on a variance blow up via the usual scaling laws that are extracted via classical log-log plots. Of course, in practical situations, one may have to deal with as many as four small, yet non-zero, parameters simultaneously: $\varepsilon$, $\sigma$, $y$, and $1/\tau$. In this context, the relative asymptotic behaviours of these parameters will matter, which is an effect already present for just pairs of small parameters~\cite{Kuehnetal}. Hence, a full mathematical analysis is beyond this work but our results show that there exists generically an open set of parameters, where inferring warning signs will be extremely difficult in practice for a class of time-correlated noise. Next, one may ask, how generic this effect of changing scaling laws is within the class of all time-correlated noise processes.}\medskip

Another natural non-Markovian option is to consider general $\alpha$-regular Volterra processes. {We denote the $\alpha$-regular Volterra processes by} $(U^{\alpha}_t)_{t\geq 0}$~\cite{Volterra}, where the parameter $\alpha\in(0,1/2)$, indicates the time correlation. The covariance of such processes is given by
\begin{align*}
\mathbb{E} [U^{\alpha}_t U^{\alpha}_s]=\int\limits_{0}^{\min\{s,t\}} K(t,r) K(s,r)~\txtd r,
\end{align*}
for a kernel $K:\mathbb{R}^2\to\mathbb{R}$ such that $| \partial_t K (t,r) | \lesssim (t-r)^{\alpha-1} \mbox{  for  } \{r<t\}$. Using the covariance function, one computes $\mathbb{E} (U^\alpha_{t} - U^\alpha_{s})^{2}  =c_{\alpha}(t-s)^{1+2\alpha}$, for an $\alpha$-dependent positive constant $c_{\alpha}$. Next, we compute a linearized fast-slow process for the variance, \textcolor{blue}{ see also~\cite{Kuehn.2020} for a similar computation in the case of a fractional Brownian motion.
In this case, regarding the scaling properties of the noise, we infer that the linearization $\xi$ along an attracting critical manifold satisfies the equation
\[
\txtd \xi_s =\frac{1}{\varepsilon} a(s) \xi_s~\txtd s +\frac{\sigma}{\varepsilon^{\alpha+\frac{1}{2}}}~\txtd U_s^\alpha,    
\]
compare~\eqref{ly:bm}.
Its solution is given by the non-autonomous Ornstein-Uhlenbeck process
\[ \xi_s =\frac{\sigma}{\varepsilon^{\alpha+\frac{1}{2}}}\int_0^s  e^{\gamma(s,u)/\varepsilon} ~\txtd U^\alpha_u,  \]
where $\gamma(s,u):=\int_u^s a(r)~\txtd r$.
Therefore, applying~\cite[Proposition 3.2]{Coupek}, we can compute the variance of $\xi$ at time $s$ as \[\text{Var}(\xi _s) =\frac{c_\alpha\sigma^2}{\varepsilon^{2\alpha+1}} \int_0^s\int_0^s e^{\gamma(s,u)/\varepsilon} e^{\gamma(s,v)/\varepsilon} |u-v|^{2\alpha-1}~\txtd u ~\txtd v .\] Finally, differentiating with respect to $s$, we obtain the nonlocal Lyapunov-type equation for the variance given by
\beann\label{lyapunov:eq}
\varepsilon V'(y) &=& 2a(y)V(y) +2 \sigma^2c_{\alpha}\cdot\\
&&\cdot\int_0^y\frac{1}{\varepsilon^{2\alpha}}
\txte^{\gamma(y,u)/\varepsilon} (y-u)^{2\alpha-1}~\txtd u.
\eeann
Using the substitution $t=\frac{y-u}{\varepsilon}$ further entails
\begin{align*}
\int\limits_{0}^{y} \frac{1}{\varepsilon^{2\alpha}} e^{\gamma(y,u)/\varepsilon} (y-u)^{2\alpha -1}~\txtd u& = \int\limits_{0}^{y/\varepsilon} e^{\gamma(y,y-\varepsilon t)/\varepsilon} t^{2\alpha -1}~\txtd t\\& \underset{\varepsilon \rightarrow 0}{\longrightarrow} \int\limits_{0}^{\infty} e^{a(y) t} t^{2\alpha -1}~\txtd t \\
& = \frac{1}{|a(y)|^{2\alpha}} \int\limits_{0}^{\infty} e^{-t} t^{2\alpha -1}~\txtd t\\
& = \frac{1}{|a(y)|^{2\alpha}} \Gamma(2\alpha),
\end{align*}
where $\Gamma(\cdot)$ is the Gamma function.}
Consequently, the limiting procedure as $\varepsilon\ra 0$ shows that the fast subsystem stationary variance in the regime before the bifurcation is given to leading-order by
\be \label{eq:varScaling_Volterra}
V_\I=\frac{\sigma^2}{|a(y)|^{2\alpha+1}}c_{\alpha}\Gamma(2\alpha).
\ee
Here, we {clearly observe} an early warning 
sign {independent of the sampling frequency} as the variance {diverges as $y\nearrow 0$. Although this scenario is different to forcing by $B_t'$, where variance tends to a finite value, the scaling law again changes in comparison to white noise forcing. In particular, for $\alpha$-regular Volterra processes} there is a whole \emph{continuum
of universal critical exponents}. The multitude of exponents can lead to practical misinterpretations for $\alpha$-regular Volterra processes. For example, if $a(y)=-2\sqrt{-y}$ as in the fold, the variance diverges as $\cO(y^{-(\alpha+1/2)})$. In combination with the assumption that there is a white noise, one may interpret the upcoming bifurcation as a pitchfork or transcritical bifurcation if $\alpha=1/2$. This can be extremely dangerous as, e.g., supercritical pitchfork bifurcations do not induce a drastic jump but a fold does. {Similar dangerous mis-interpretations can evidently arise if one aims to predict the distance to the bifurcation point via scaling laws if the precise type of noise present in data is unknown.} 

There are several important examples for $\alpha$-regular Volterra processes. The most famous one is fractional Brownian motion 
\begin{align*}
	    B^H_t = c(\cH) 
	    \int\limits_{\mathbb{R}} \int\limits_{0}^{t} (r-s)^{\cH-3/2}_{+}~\txtd r~\txtd W_{s}
\end{align*}
with Hurst index $\cH\in(1/2,1)$, where $c(\cH)$ is a $\cH$-dependent positive constant, which ensures that $\mathbb{E} (B^{\cH}_{1})^2=1$. The previous computation holds for $\cH:=\alpha+1/2$ (see also \cite{Kuehn.2020}). Another example is the Rosenblatt process (see \cite{Tudor.2008}), which is a non-Gaussian processes with stationary increments and the same covariance as fractional Brownian motion
\begin{align*}
	R^{H}_{t} = \frac{\sqrt{\cH}(2\cH-1)}{\beta(\cH/2,1-\cH)}\cdot\\
	 \int\limits_{\mathbb{R}^{2}} \Bigg( \int\limits_{0}^{t} (u-y_{1})^{-\frac{2-\cH}{2}}_{+} (u-y_{2})_{+}^{-\frac{2-\cH}{2}}~\txtd u \Bigg)~\txtd W_{y_{1}}~\txtd W_{y_{2}},
\end{align*}
where $\beta$ stands for the beta function. Fractional Brownian motion and the Rosenblatt process are both self-similar with parameter $\cH\in(1/2,1)$. Their stationary variances in the regime before the bifurcation can be inferred from \eqref{eq:varScaling_Volterra} by using the relation $\cH=\alpha +\nicefrac{1}{2}$ (see~\cite{Maslowski.2017}) as both processes are of $\alpha$-regular Volterra type, so
\begin{align}
	V_\I = \frac{\sigma^2 \cH\Gamma(2\cH)}{a(y)^{2\cH}}. \label{eq:scalingLaw_fBM}
\end{align}

Using the linearizations $a(y)$ in \eqref{eq:a_pf_tc_lp}, different scaling laws arise for pitchfork, transcritical and fold bifurcations. Via formulas \eqref{eq:scalingLaw_color} and $\eqref{eq:scalingLaw_fBM}$, we get Table \ref{tab:scaling_noiseTypes}.
\begin{table}[h]
	\centering
	\begin{tabular}{l|c|c|c}
		& pitchfork & transcritical & fold \\
		\hline
	(C1) white noise	& -1 & -1 & -\nicefrac{1}{2} \\
		\hline
	(C2) coloured noise	& 0 & 0 & 0 \\
		\hline
	(C3) fractional Brownian motion	& -2$\cH$ & -2$\cH$ & -$\cH$ \\
		\hline
	(C4) Rosenblatt process	& -2$\cH$ & -2$\cH$ & -$\cH$ \\
\end{tabular}
\caption{Scaling laws for classical bifurcation types for the four different noise types.}
\label{tab:scaling_noiseTypes}
\end{table}

Despite their structural difference, fractional Brownian motion and the Rosenblatt process obey the same variance scaling near the bifurcation point as they share the same covariance function.\medskip

\textbf{Climate Tipping Numerics.} To illustrate the impact of our results, we consider a stochastic Stommel-Cessi \cite{Cessi.1994} box model for the Atlantic Meridional Overturning Circulation (AMOC). It is representative in the sense that its S-shaped fold bifurcation with corresponding tipping points appears also in the context of global energy balance models \cite{Ghil.2015,Ghil.2020} in the climate literature. Although being very basic, box models are powerful tools some of which can capture several key dynamics and ocean
properties \cite{Alkhayuon.2019,Wood.2019}. Here, we focus on the stochastic Stommel-Cessi box model, which is a two-dimensional fast-slow system given by (see \cite[Sec.\ 6.2.1]{Gentz.2006} and \cite{KuehnCT2})
\begin{subequations} \label{eq:Stommel_timeCorrNoise}
	\begin{align} 
		x' &= y - x\left(1+\eta^2\cdot(1-x)^2\right) + \sigma~C'_t, \label{eq:Stommel_fast}\\
		y' &= -\varepsilon \label{eq:Stommel_slow}
	\end{align}
\end{subequations}
where $x$ represents the salinity difference.

\begin{figure}[h]
	\centering
	\subfloat[\ Sample paths along critical manifold for white noise{: the movement of $100$ sample paths along the upper part of the critical manifold is shown before they jump off to the lower attracting branch of the critical manifold. The zoomed in version shows the paths structure resulting from white noise forcing.} \label{fig:pathsCMa_Stommel_BM_T45yIni1K1643pIni1K4sigma0K01eps-0K01}]{\includegraphics[width=0.5\textwidth]{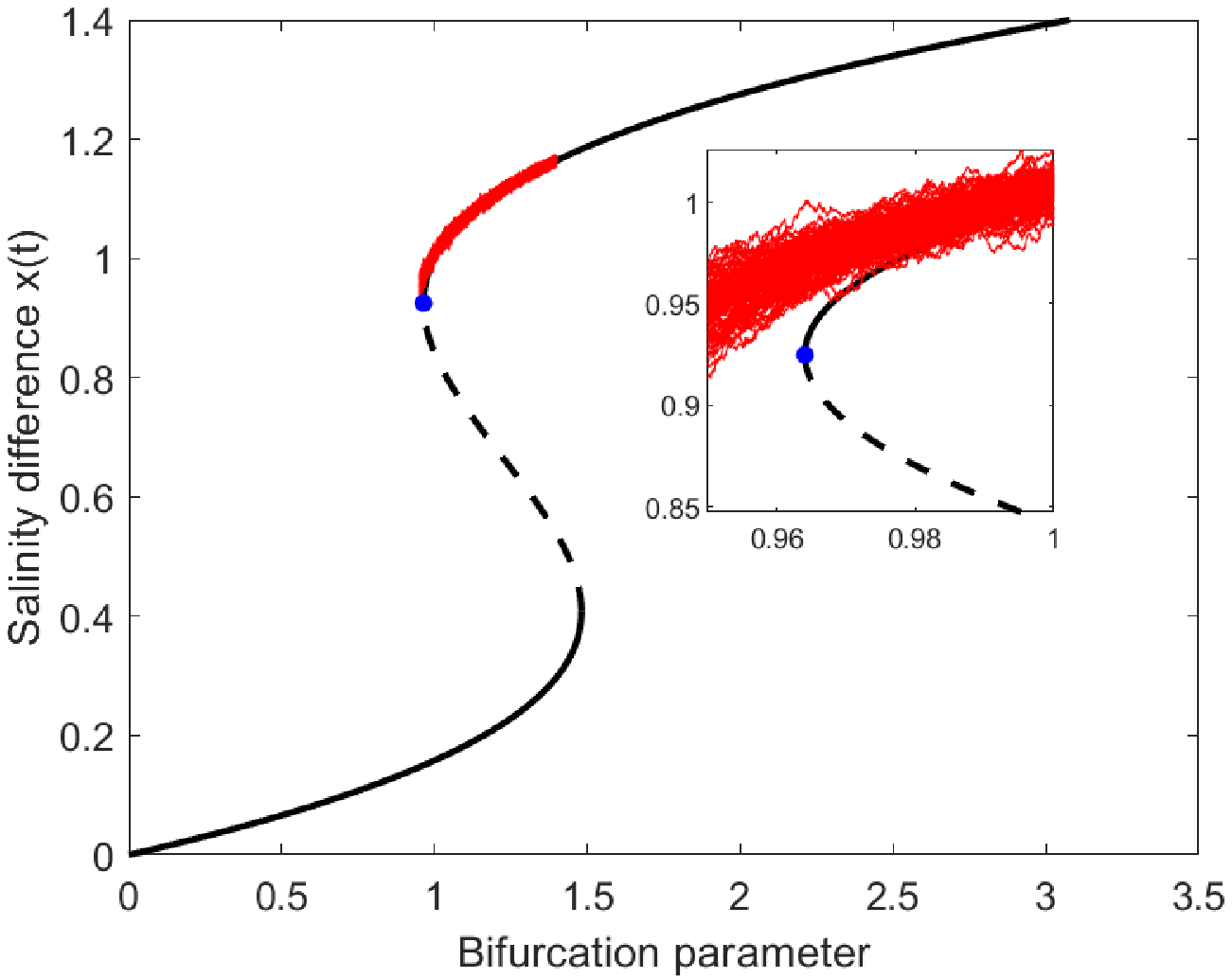}}\hfill
	\subfloat[\ Log-log-plot of {the} variance {evolution (in blue) of $x(t)$ with} decreasing distance to bifurcation{: sample size $M=10^4$, time discretization $\Delta t=10^{-3}$, final time $T=45$. The red straight line is a linear best fit in a least-squares sense via the MATLAB function \textit{polyfit}\protect\footnotemark.} \label{fig:loglog_var_MC_Stommel_BM_T45yIni1K1643pIni1K4sigma0K01eps-0K01}]{\includegraphics[width=0.45\textwidth]{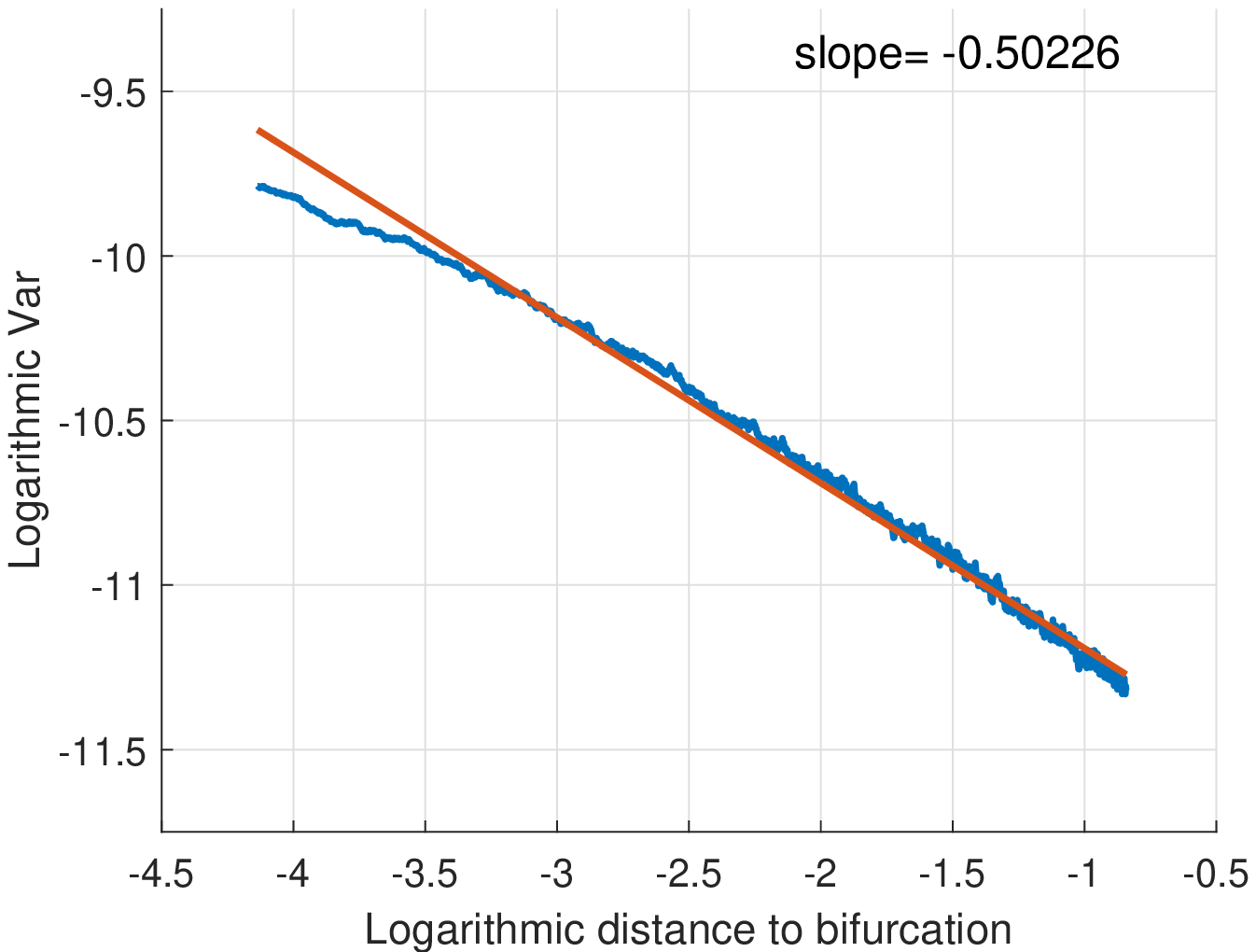}}
	\caption{(C1) White noise }
	\label{fig:Stommel_BM_T45yIni1K1643pIni1K4sigma0K01eps-0K01}
\end{figure}
\footnotetext{{source: \url{https://de.mathworks.com/help/matlab/ref/polyfit.html}}}

\begin{figure}[h]
	\centering
	\subfloat[\ Sample paths along critical manifold for coloured noise{: the movement of $100$ sample paths along the upper part of the critical manifold is shown before they jump off to the lower attracting branch of the critical manifold. The zoomed in version shows the paths structure resulting from coloured noise forcing.} \label{fig:pathsCMa_Stommel_color_tau0K05_T45yIni1K1643pIni1K4sigma0K01eps-0K01}]{\includegraphics[width=0.5\textwidth]{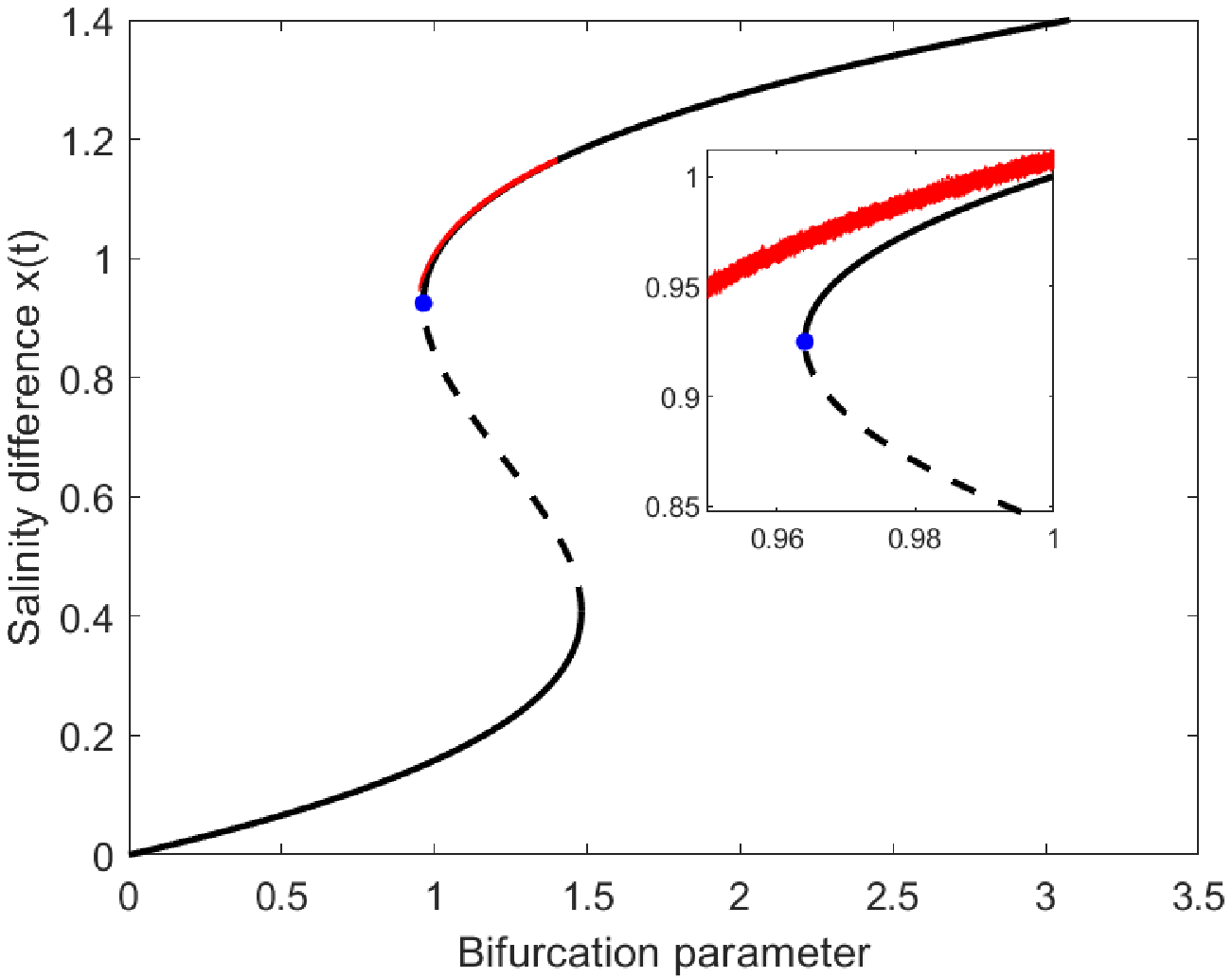}}\hfill
	\subfloat[\  Log-log-plot of {the} variance {evolution (in blue) of $x(t)$ with} decreasing distance to bifurcation{: sample size $M=10^4$, time discretization $\Delta t=10^{-3}$, final time $T=45$. The red straight line is a linear best fit in a least-squares sense via the MATLAB function \textit{polyfit}.} \label{fig:loglog_var_MC_Stommel_color_tau0K05_T45yIni1K1643pIni1K4sigma0K01eps-0K01}]{\includegraphics[width=0.45\textwidth]{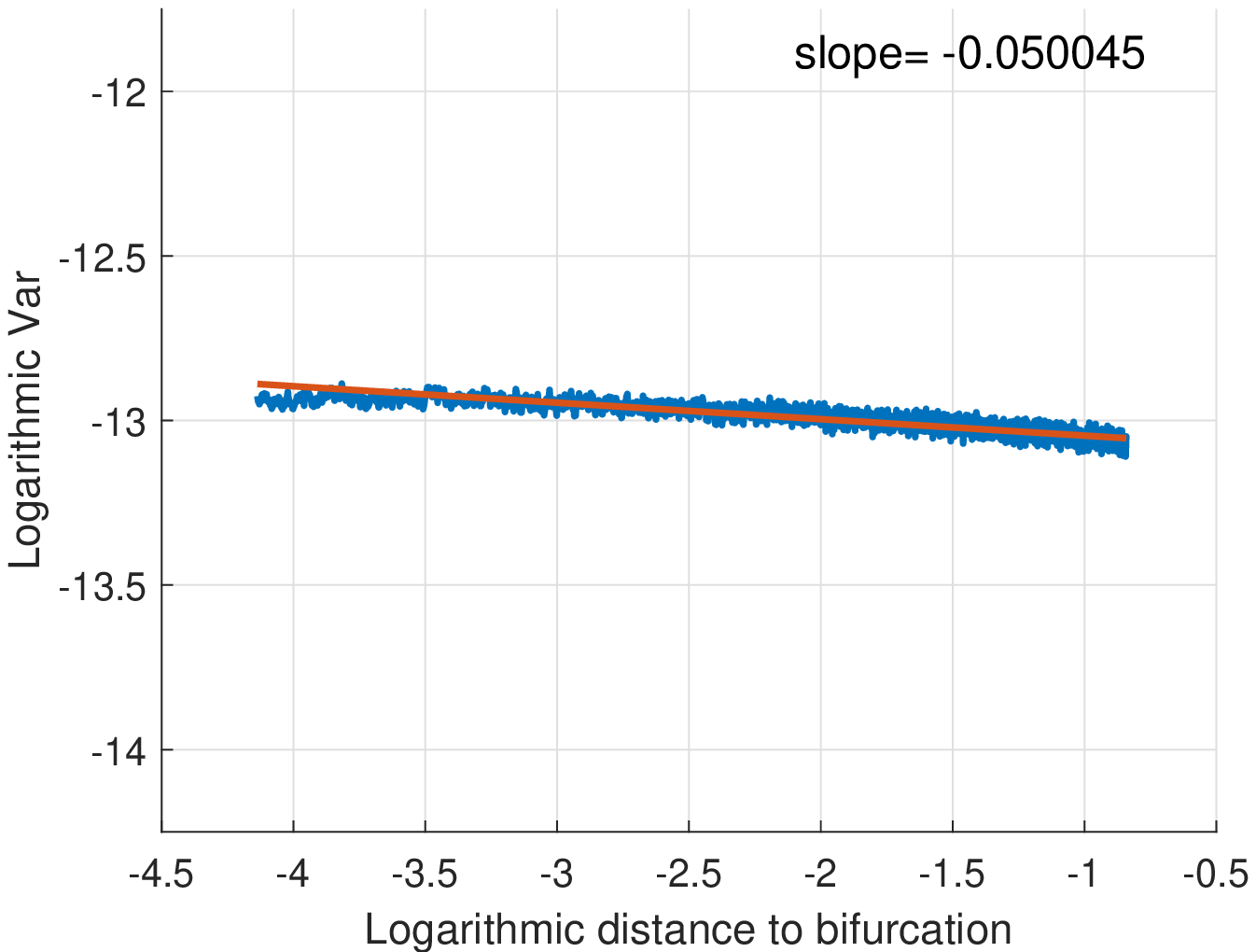}}
	\caption{(C2) Coloured noise with correlation $\tau=0.05$. }
	\label{fig:Stommel_color_T45yIni1K1643pIni1K4sigma0K01eps-0K01}
\end{figure}

\begin{figure}[h]
	\centering
	\subfloat[\ Sample paths along critical manifold for fractional BM{: the movement of $100$ sample paths along the upper part of the critical manifold is shown before they jump off to the lower attracting branch of the critical manifold. The zoomed in version shows the paths structure resulting from forcing via a fractional BM.} \label{fig:pathsCMa_Stommel_fBM_H0K9_T45yIni1K1643pIni1K4sigma0K01eps-0K01}]{\includegraphics[width=0.5\textwidth]{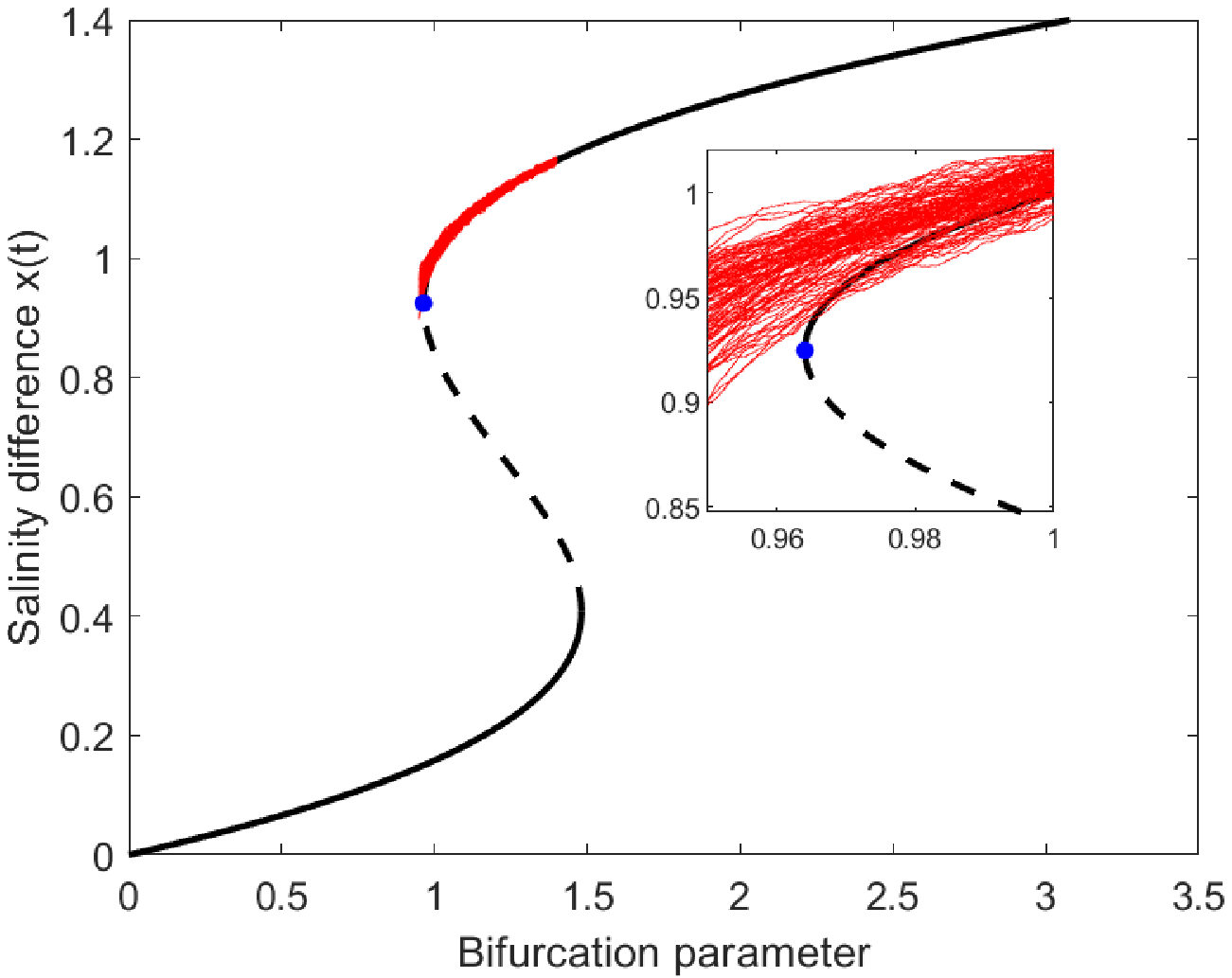}}\hfill
	\subfloat[\  Log-log-plot of {the} variance {evolution (in blue) of $x(t)$ with} decreasing distance to bifurcation{: sample size $M=10^4$, time discretization $\Delta t=10^{-3}$, final time $T=45$. The red straight line is a linear best fit in a least-squares sense via the MATLAB function \textit{polyfit}.} \label{fig:loglog_var_MC_Stommel_fBM_H0K9_T45yIni1K1643pIni1K4sigma0K01eps-0K01}]{\includegraphics[width=0.45\textwidth]{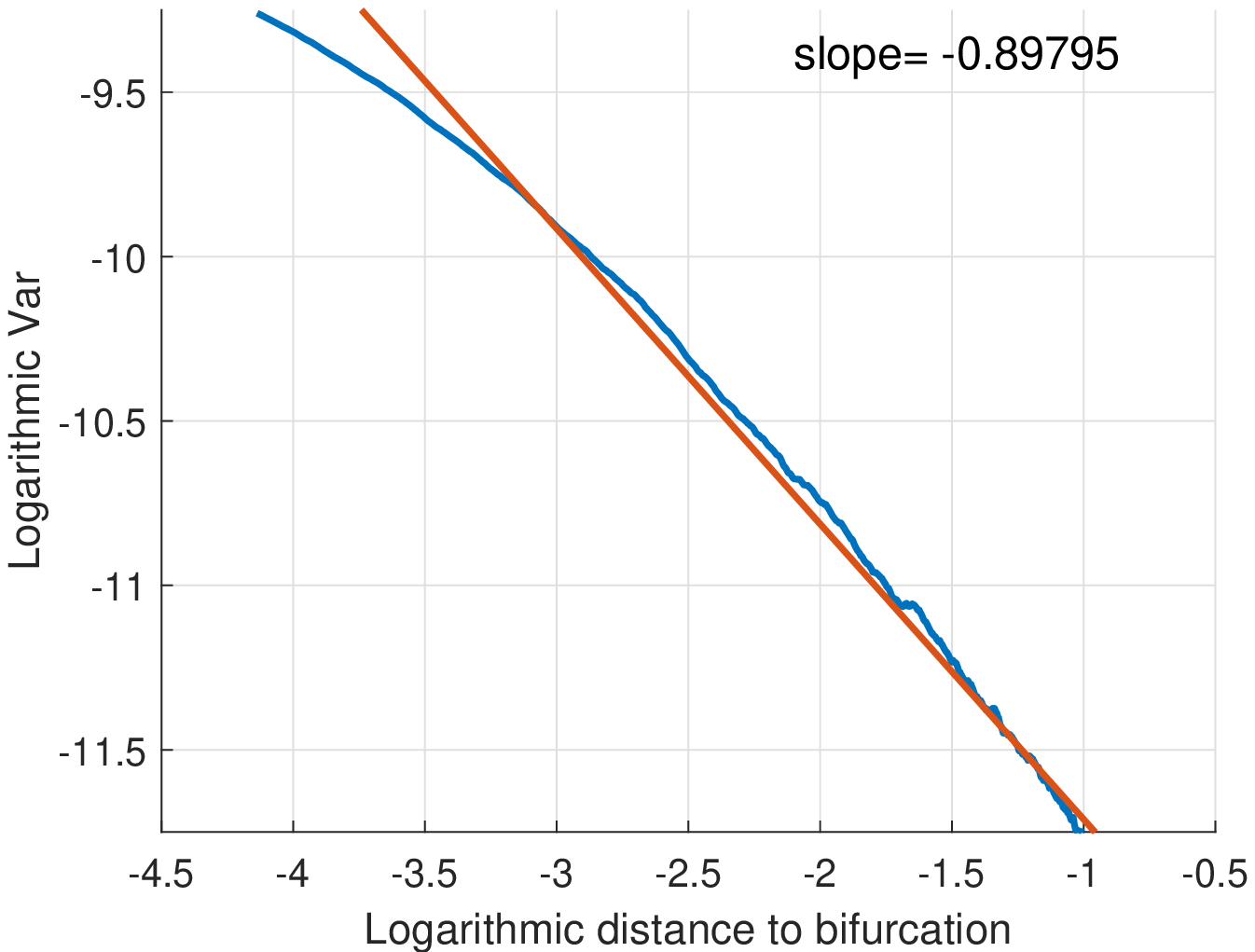}}
	\caption{(C3) Fractional Brownian Motion with $\cH=0.9$. }
	\label{fig:Stommel_fBM_T45yIni1K1643pIni1K4sigma0K01eps-0K01}
\end{figure}

To model time correlation in the stochastic forcing, we take $\left(C_t\right)_{t \in [0,T]}$ as one of the processes from Table \ref{tab:scaling_noiseTypes}. The ratio between diffusive and advective time scale is fixed to $\eta^2=7.5$ {(value taken from \cite{Cessi.1994})}.
The dynamic bifurcation parameter $y$ is proportional to the atmospheric freshwater flux. The deterministic critical manifold has a classical $S$-shaped form
\begin{align} \label{eq:critManifold_Stommel}
	C_0 &= \{(x,y) \in \R^2: y = x(1+7.5(1-x)^2) \}.
\end{align}
with two fast-subsystem fold bifurcations. We are interested in the fold point in $P_{\text{fold}} = \left(\nicefrac{1}{15}(10+\sqrt{15}),\nicefrac{11}{9}-\nicefrac{1}{\sqrt{15}}\right)$. {We use the computing platform MATLAB}\footnotetext{{source: https://de.mathworks.com/products/matlab.html}}{, version 2019a, on a Desktop PC (64 Bit-operating system Windows 7, Intel\textregistered \ Core\texttrademark \ i7-3770 CPU with 3.4 GHz, 16GB RAM).} We fix $\varepsilon=0.01${, $\sigma=0.01$}, and $(x(0),y(0)) = (x_0,1.4)$ on the attracting upper branch $C_0^{a,+}$. States on $C_0^{a,+}$ represent {the} AMOC {in a weak state}. We sample $M=10^4$ paths of $\left(x(t)\right)_{t\in [0,T]}$ of \eqref{eq:Stommel_timeCorrNoise} numerically and calculate the sample variance $\var(x(t))$ based on these time series up to $T=45$. {Note that for practical climate data, where often only one sample path is available, the adaptation of a classical sliding window approach for the variance calculation might be a remedy. However, this is not straightforward in our case due to the inherent time correlation showing yet another practical difficulty, when noise is time-correlated.}

\begin{figure}[h]
	\centering
	\subfloat[\ Sample paths along critical manifold for Rosenblatt process:{: the movement of $100$ sample paths starting in $y_0=1.0642$ along the upper part of the critical manifold is shown before they jump off to the lower attracting branch of the critical manifold. The zoomed in version shows the paths structure resulting from forcing via a Rosenblatt process.} \label{fig:pathsCMa_Stommel_scaledRoPro_H0K9_T10yIni1K0469pIni1K0642sigma0K01eps-0K01dt0K01}]{\includegraphics[width=0.5\textwidth]{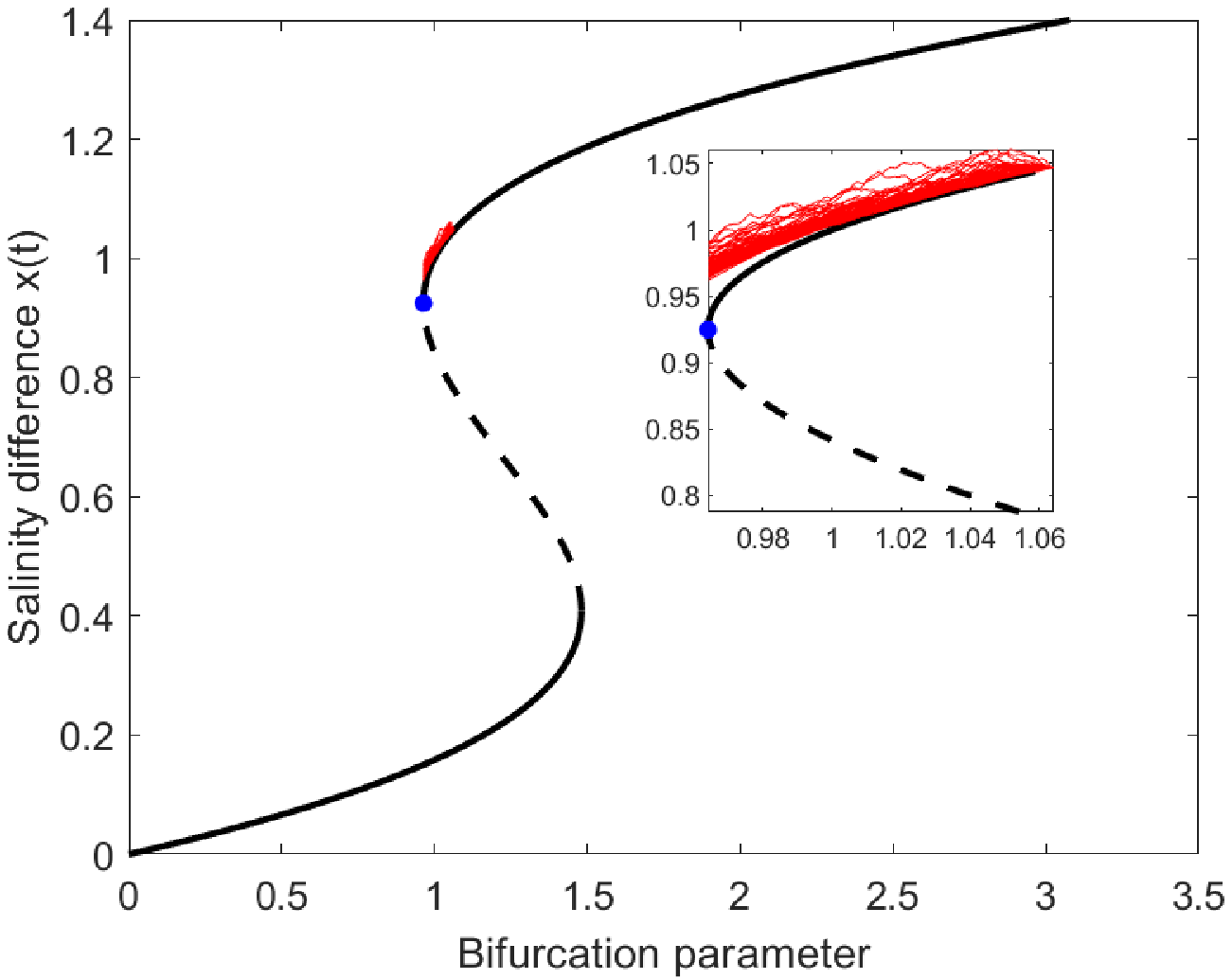}}\hfill
	\subfloat[\  Log-log-plot of {the} variance {evolution (in blue) of $x(t)$ with} decreasing distance to bifurcation{: sample size $M=10^3$, time discretization $\Delta t=10^{-2}$, final time $T=10$. The red straight line is a linear best fit in a least-squares sense via the MATLAB function \textit{polyfit}. Note that the simulation of the Rosenblatt process is more involved and that a coarser time discretization and fewer sample paths have been used resulting in a poorer performance of the linear fit.} \label{fig:loglog_var_MC_Stommel_sclaedRoPro_H0K9_T10yIni1K0469pIni1K0642sigma0K01eps-0K01dt0K01}]{\includegraphics[width=0.45\textwidth]{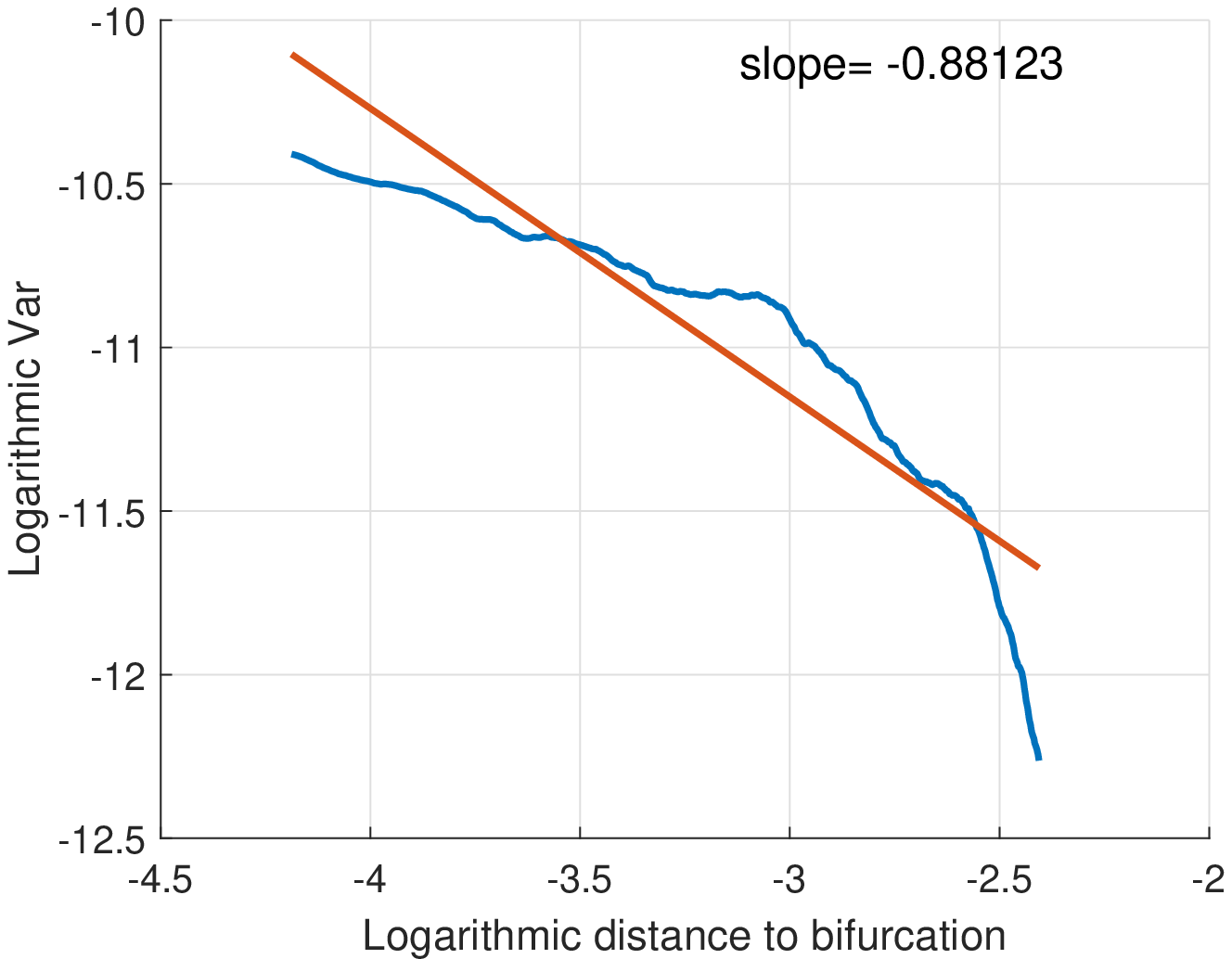}}
	\caption{(C4) Rosenblatt process with $\cH=0.9$. }
	\label{fig:Stommel_scaledRoPro_H0K9_T10yIni1K0469pIni1K0642sigma0K01eps-0K01dt0K01}
\end{figure}

Although standard techniques apply for the simulation of the noise types (C1)-(C3)~\cite{Asmussen.2007,Kloeden.1992}, the simulation of the Rosenblatt process turns out to be challenging. We use a discretization of the finite-time stochastic integral representation \cite[Prop.\ 1]{Tudor.2008}. Limitations with respect to memory allocation might arise. Therefore, we restrict our simulation to $T=10$, $\Delta t=10^{-2}$, $M=10^3$, and use the initial condition $(x_0,y_0)=(x_0,1.0642)$ on $C_0^{a,+}$. The precise quantitative variance computation of the Rosenblatt process is a numerically challenging task as already pointed out in \cite{Tudor.2010} but the scaling law can still be recovered. In Figures \ref{fig:Stommel_BM_T45yIni1K1643pIni1K4sigma0K01eps-0K01}-\ref{fig:Stommel_scaledRoPro_H0K9_T10yIni1K0469pIni1K0642sigma0K01eps-0K01dt0K01}, the evolution of the sample paths along the critical manifold and the numerical rates of the variance scalings in $\log$-$\log$-scale are depicted for all four noise types. We observe a close agreement of the numerical rates with the theoretically proven ones in Table \ref{tab:scaling_noiseTypes}. {Since there are different scaling laws possible for the different time-correlated noises, ranging from a practically extremely difficult to detect variance change to a clear growth in variance, predictability becomes very difficult if the noise type in data is not known. Therefore, it is important to extend the debates from the literature regarding white noise and EWS for climate tipping events~\cite{Dakosetal,DitlevsenJohnsen}, to time-correlated processes. These provide another possible explanation for discrepancies in predictability.}.

The code used for generating Figures \ref{fig:Stommel_BM_T45yIni1K1643pIni1K4sigma0K01eps-0K01}-\ref{fig:Stommel_scaledRoPro_H0K9_T10yIni1K0469pIni1K0642sigma0K01eps-0K01dt0K01} can be found at \url{https://github.com/kerstinLux/EWStimeCorrelatedNoise}.

\bigskip

\textbf{Conclusions and Outlook.}
To our best knowledge, this is the first work which systematically investigates EWS for time-correlated non-Markovian processes providing a precise scaling dependence of the warning signs upon the time-correlation and self-similarity of the noise. We analysed the influence of the noise for the prediction of bifurcation points and found that coloured noise can make us ``blind'' for tipping in the sense that {no divergence of the variance}
can be observed before passing the tipping point. Also for other types of non-Markovian noise such as Volterra-type processes, and specifically fractional Brownian motion, the warning sign scaling laws change very significantly. We cross-validated our theoretical findings numerically for a climate box model.
It would be interesting to analyse the effect of time-correlated noise on more realistic models of the AMOC such as box models allowing for a parameter calibration according to general circulation models \cite{Alkhayuon.2019,Wood.2019}. Furthermore, it would be relevant to systematically re-evaluate existing time series analysis results regarding warning signs in cases, where there seems to be ongoing debates in the applied literature, when and if warning signs are present~\cite{BoettingerRossHastings,Dakosetal,DitlevsenJohnsen}. Taking the viewpoint of allowing for memory effects could yield significant additional insights to existing arguments.\\

\textbf{Acknowledgments:} CK and KL acknowledge support via the TiPES project funded by the European Unions Horizon 2020 research and innovation programme under grant agreement No. 820970. This is TiPES publication \#119. CK acknowledges partial support of the VolkswagenStiftung via a Lichtenberg Professorship.

\bibliographystyle{siam}
\bibliography{myLiterature}
\end{document}